\documentclass[12pt]{article}

\usepackage{amsmath,amssymb,amsthm,amscd,a4wide}

\begin{document}

\newcommand{\End}{{\rm{End}\ts}}
\newcommand{\Hom}{{\rm{Hom}}}
\newcommand{\Mat}{{\rm{Mat}}}
\newcommand{\ch}{{\rm{ch}\ts}}
\newcommand{\chara}{{\rm{char}\ts}}
\newcommand{\diag}{ {\rm diag}}
\newcommand{\non}{\nonumber}
\newcommand{\wt}{\widetilde}
\newcommand{\wh}{\widehat}
\newcommand{\ot}{\otimes}
\newcommand{\la}{\lambda}
\newcommand{\La}{\Lambda}
\newcommand{\De}{\Delta}
\newcommand{\al}{\alpha}
\newcommand{\be}{\beta}
\newcommand{\ga}{\gamma}
\newcommand{\Ga}{\Gamma}
\newcommand{\ep}{\epsilon}
\newcommand{\ka}{\kappa}
\newcommand{\vk}{\varkappa}
\newcommand{\si}{\sigma}
\newcommand{\vp}{\varphi}
\newcommand{\de}{\delta}
\newcommand{\ze}{\zeta}
\newcommand{\om}{\omega}
\newcommand{\ee}{\epsilon^{}}
\newcommand{\su}{s^{}}
\newcommand{\hra}{\hookrightarrow}
\newcommand{\ve}{\varepsilon}
\newcommand{\ts}{\,}
\newcommand{\vac}{\mathbf{1}}
\newcommand{\di}{\partial}
\newcommand{\qin}{q^{-1}}
\newcommand{\tss}{\hspace{1pt}}
\newcommand{\Sr}{ {\rm S}}
\newcommand{\U}{ {\rm U}}
\newcommand{\BL}{ {\overline L}}
\newcommand{\BE}{ {\overline E}}
\newcommand{\BP}{ {\overline P}}
\newcommand{\AAb}{\mathbb{A}\tss}
\newcommand{\CC}{\mathbb{C}\tss}
\newcommand{\KK}{\mathbb{K}\tss}
\newcommand{\QQ}{\mathbb{Q}\tss}
\newcommand{\SSb}{\mathbb{S}\tss}
\newcommand{\ZZ}{\mathbb{Z}\tss}
\newcommand{\X}{ {\rm X}}
\newcommand{\Y}{ {\rm Y}}
\newcommand{\Z}{{\rm Z}}
\newcommand{\Ac}{\mathcal{A}}
\newcommand{\Lc}{\mathcal{L}}
\newcommand{\Mc}{\mathcal{M}}
\newcommand{\Pc}{\mathcal{P}}
\newcommand{\Qc}{\mathcal{Q}}
\newcommand{\Tc}{\mathcal{T}}
\newcommand{\Sc}{\mathcal{S}}
\newcommand{\Bc}{\mathcal{B}}
\newcommand{\Ec}{\mathcal{E}}
\newcommand{\Fc}{\mathcal{F}}
\newcommand{\Hc}{\mathcal{H}}
\newcommand{\Uc}{\mathcal{U}}
\newcommand{\Vc}{\mathcal{V}}
\newcommand{\Wc}{\mathcal{W}}
\newcommand{\Ar}{{\rm A}}
\newcommand{\Br}{{\rm B}}
\newcommand{\Ir}{{\rm I}}
\newcommand{\Fr}{{\rm F}}
\newcommand{\Jr}{{\rm J}}
\newcommand{\Or}{{\rm O}}
\newcommand{\GL}{{\rm GL}}
\newcommand{\Spr}{{\rm Sp}}
\newcommand{\Rr}{{\rm R}}
\newcommand{\Zr}{{\rm Z}}
\newcommand{\gl}{\mathfrak{gl}}
\newcommand{\middd}{{\rm mid}}
\newcommand{\Pf}{{\rm Pf}}
\newcommand{\Norm}{{\rm Norm\tss}}
\newcommand{\oa}{\mathfrak{o}}
\newcommand{\spa}{\mathfrak{sp}}
\newcommand{\osp}{\mathfrak{osp}}
\newcommand{\g}{\mathfrak{g}}
\newcommand{\h}{\mathfrak h}
\newcommand{\n}{\mathfrak n}
\newcommand{\z}{\mathfrak{z}}
\newcommand{\Zgot}{\mathfrak{Z}}
\newcommand{\p}{\mathfrak{p}}
\newcommand{\sll}{\mathfrak{sl}}
\newcommand{\agot}{\mathfrak{a}}
\newcommand{\qdet}{ {\rm qdet}\ts}
\newcommand{\Ber}{ {\rm Ber}\ts}
\newcommand{\HC}{ {\mathcal HC}}
\newcommand{\cdet}{ {\rm cdet}}
\newcommand{\tr}{ {\rm tr}}
\newcommand{\str}{ {\rm str}}
\newcommand{\loc}{{\rm loc}}
\newcommand{\Gr}{{\rm G}}
\newcommand{\sgn}{ {\rm sgn}\ts}
\newcommand{\ba}{\bar{a}}
\newcommand{\bb}{\bar{b}}
\newcommand{\bi}{\bar{\imath}}
\newcommand{\bj}{\bar{\jmath}}
\newcommand{\bk}{\bar{k}}
\newcommand{\bl}{\bar{l}}
\newcommand{\hb}{\mathbf{h}}
\newcommand{\Sym}{\mathfrak S}
\newcommand{\fand}{\quad\text{and}\quad}
\newcommand{\Fand}{\qquad\text{and}\qquad}
\newcommand{\For}{\qquad\text{or}\qquad}
\newcommand{\OR}{\qquad\text{or}\qquad}

\renewcommand{\theequation}{\arabic{section}.\arabic{equation}}

\newtheorem{thm}{Theorem}[section]
\newtheorem{lem}[thm]{Lemma}
\newtheorem{prop}[thm]{Proposition}
\newtheorem{cor}[thm]{Corollary}
\newtheorem{conj}[thm]{Conjecture}
\newtheorem*{mthma}{Theorem A}
\newtheorem*{mthmb}{Theorem B}

\theoremstyle{definition}
\newtheorem{defin}[thm]{Definition}

\theoremstyle{remark}
\newtheorem{remark}[thm]{Remark}
\newtheorem{example}[thm]{Example}

\newcommand{\bth}{\begin{thm}}
\renewcommand{\eth}{\end{thm}}
\newcommand{\bpr}{\begin{prop}}
\newcommand{\epr}{\end{prop}}
\newcommand{\ble}{\begin{lem}}
\newcommand{\ele}{\end{lem}}
\newcommand{\bco}{\begin{cor}}
\newcommand{\eco}{\end{cor}}
\newcommand{\bde}{\begin{defin}}
\newcommand{\ede}{\end{defin}}
\newcommand{\bex}{\begin{example}}
\newcommand{\eex}{\end{example}}
\newcommand{\bre}{\begin{remark}}
\newcommand{\ere}{\end{remark}}
\newcommand{\bcj}{\begin{conj}}
\newcommand{\ecj}{\end{conj}}

\newcommand{\bal}{\begin{aligned}}
\newcommand{\eal}{\end{aligned}}
\newcommand{\beq}{\begin{equation}}
\newcommand{\eeq}{\end{equation}}
\newcommand{\ben}{\begin{equation*}}
\newcommand{\een}{\end{equation*}}

\newcommand{\bpf}{\begin{proof}}
\newcommand{\epf}{\end{proof}}

\def\beql#1{\begin{equation}\label{#1}}

\title{\Large\bf Casimir elements from the Brauer--Schur--Weyl duality}

\author{{N. Iorgov,\quad A. I. Molev\quad and\quad
E. Ragoucy}}

\date{} % Start April 2012
\maketitle

\vspace{25 mm}

\begin{abstract}
We consider Casimir elements for the orthogonal and symplectic
Lie algebras constructed with the use of the Brauer algebra.
We calculate the images of these elements under the Harish-Chandra
isomorphism and thus show that they (together with the Pfaffian-type element
in the even orthogonal case) are algebraically independent
generators of the centers
of the corresponding universal enveloping algebras.

\vspace{5 pt}

Preprint LAPTH-027/12
\end{abstract}

%%%\vspace{5 mm}
%%%
%%%{\it Key words:}
%%%

\vspace{30 mm}

\noindent
Bogolyubov Institute for Theoretical Physics\newline
03680, Kyiv, Ukraine\newline
iorgov@bitp.kiev.ua

\vspace{7 mm}

\noindent
School of Mathematics and Statistics\newline
University of Sydney,
NSW 2006, Australia\newline
alexander.molev@sydney.edu.au

\vspace{7 mm}

\noindent
LAPTH, Chemin de Bellevue, BP 110\newline
CNRS and Universit\'e de Savoie\newline
F-74941 Annecy-le-Vieux cedex, France\newline
ragoucy@lapp.in2p3.fr

\newpage

\section{Introduction}
\label{sec:int}
\setcounter{equation}{0}

It is well-known that the Schur--Weyl duality can be used to get natural
constructions of families of Casimir elements for
the classical Lie algebras. For some particular choices of parameters
the images of such elements under the Harish-Chandra isomorphism can
be calculated in an explicit form. For the general linear Lie algebras $\gl_N$
this leads to an explicit construction of a linear basis of the center
of the universal enveloping algebra $\U(\gl_N)$. The basis elements are known
as the quantum immanants and their Harish-Chandra images are the factorial (or shifted)
Schur functions; see \cite{o:qi} and \cite{oo:ss}. Constructing quantum
immanant-type bases of the centers of the universal enveloping algebras $\U(\oa_N)$
and $\U(\spa_N)$ for the orthogonal and symplectic Lie algebras remains
an open problem; see, however \cite{mr:cm}, \cite{n:ce},
\cite{oo:ss2} and \cite{o:gs} for some results in that direction.

In this paper we consider
generators of the centers of $\U(\oa_N)$ and $\U(\spa_N)$
obtained by an application of the
Brauer--Schur--Weyl duality. They are associated with one-dimensional
representations of the Brauer algebra and take the form of some versions
of noncommutative determinants and permanents.
We give explicit formulas for the
Harish-Chandra images of these elements; the images turn out to coincide
with the factorial (or double) complete and elementary
symmetric functions.

In more detail, we regard $\CC^N$ as the vector representation
of each of the groups $\Or_N$ and $\Spr_N$ (the latter with even $N$).
The space of tensors
\beql{tenprk}
\underbrace{\CC^{N}\ot\dots\ot\CC^{N}}_m
\eeq
carries the diagonal action of each group. By the Schur--Weyl
duality, the centralizer of the action
of the orthogonal group $\Or_N$ or symplectic group $\Spr_N$ in $\End(\CC^N)^{\ot m}$
is generated by the homomorphic image
of the action of the respective Brauer algebra $\Bc_m(N)$ or $\Bc_m(-N)$.
Consider multiple tensor products
\beql{tenprku}
\U(\g_N)\ot\underbrace{\End\CC^{N}\ot\dots\ot\End\CC^{N}}_m\ts,
\eeq
where $\g_N$ denotes one of the Lie algebras $\oa_N$ or $\spa_N$.
We let $F=[F_{ij}]$ denote the $N\times N$ matrix
whose entries $F_{ij}$ are the standard generators of $\g_N$; see
definitions in Sec.~\ref{sec:cas}.
For each $a=1,\dots,m$ we denote by $F_a$ the respective element of the algebra
\eqref{tenprku},
\beql{ea}
F_a=\sum_{i,j=1}^N F_{ij}\ot 1^{\ot(a-1)}\ot e_{ij}\ot 1^{\ot(m-a)},
\eeq
where $e_{ij}\in\End\CC^{N}$ denote the standard matrix units.
We regard an arbitrary element $C$ of the respective
Brauer algebra $\Bc_m(N)$ or $\Bc_m(-N)$ as
an operator in the space \eqref{tenprk}.
The adjoint action of the respective group $\Gr_N=\Or_N$ or $\Gr_N=\Spr_N$
on the corresponding Lie algebra $\g_N$
amounts to conjugations of the matrix $F$ by elements of the group.
Hence, since the action of the group $\Gr_N$ on the space \eqref{tenprk}
commutes with the action
of the Brauer algebra,
we find that the elements
\beql{tracegn}
\tr\ts C\tss (F_1+u_1)\dots (F_m+u_m),
\eeq
with the trace taken over all $m$ copies of $\End\CC^{N}$, belong
to the subalgebra $\U(\g_N)^{\Gr_N}$ of $\Gr_N$-invariants in $\U(\g_N)$.
We are using the notation $F+u$ to
indicate the matrix $F+u\tss1$.
If $\Gr_N=\Or_N$ with odd $N$ or $\Gr_N=\Spr_N$ with even $N$, then
the subalgebra $\U(\g_N)^{\Gr_N}$
coincides with the center of $\U(\g_N)$. If $\Gr_N=\Or_N$ with even $N$,
then $\U(\g_N)^{\Gr_N}$ is a proper subalgebra of the center which contains
an additional Pfaffian-type Casimir element.

Note that the quantum immanants of \cite{o:qi}
are elements of the form \eqref{tracegn},
where $F$ should be replaced with the matrix $E=[E_{ij}]$ formed by the basis
elements $E_{ij}$ of the general linear Lie algebra,
$C$ is a primitive idempotent of $\CC[\Sym_m]$ associated with
a standard tableau and the $u_i$ are contents of the tableau.

In a recent work \cite{m:ff} an explicit construction of
generators of the center of the affine vertex algebra $V(\g_N)$
at the critical level was given. Here $V(\g_N)$
is the vacuum module over the affine Kac--Moody algebra $\wh\g_N$.
Under the evaluation homomorphism
the generators of the center of $V(\g_N)$ get mapped into Casimir elements
of the form \eqref{tracegn}. More precisely,
in the orthogonal case
the images of the central elements defined in \cite{m:ff} have the form
of a differential operator
\beql{tracedi}
\tr\ts S^{(m)}\tss (-\di_t+F_1\tss t^{-1})\dots (-\di_t+F_m\tss t^{-1}),
\eeq
whose coefficients are Casimir elements, where $S^{(m)}$ denotes the
symmetrizer in the Brauer algebra. In the symplectic case with $N=2n$
the images of the central elements
have the form of \eqref{tracedi} with the additional factor $(n-m+1)^{-1}$
and the values of $m$ restricted to $1\leqslant m\leqslant 2n$.
Multiplying the operator \eqref{tracedi}
from the left by $t^m$ we get an expression of the form \eqref{tracegn}
with $C=S^{(m)}$ and the parameters specialized as $u_i=u+m-i$
for $i=1,\dots,m$, where $u=-t\tss\di_t$.

To identify this family of central elements
we calculate their eigenvalues in highest weight representations of $\g_N$
which is equivalent to
finding their images under the Harish-Chandra isomorphism.
The key starting point in the orthogonal case
$\g_N=\oa_N$ with $N=2n$ or $N=2n+1$ is Theorem~\ref{thm:caseven}
which implies that the Casimir
elements
\beql{genb}
\tr\ts S^{(2k)}\tss (F_1+k-1)\dots (F_{2k}-k),\qquad k=1,\dots,n
\eeq
are algebraically independent generators of the algebra of invariants
$\U(\oa_N)^{\Or_N}$. Their eigenvalues in the highest weight
representations with the highest weight $(\la_1,\dots,\la_n)$
coincide with the factorial
complete symmetric functions $h_k(l_1^2,\dots,l_n^2\ts|\ts a)$,
where the $l_i$ are the labels of the highest weight
shifted by the half-sum of the positive roots and
$a=(a_i)$ is a sequence of parameters; see Sec.~\ref{subsec:cps}.
In the symplectic case $\g_N=\spa_{2n}$ the Casimir
elements defined by \eqref{genb} with a normalization factor
are algebraically independent generators of the center of $\U(\spa_{2n})$.
Their eigenvalues in the highest weight
representations
coincide with the factorial
elementary symmetric functions $e_k(l_1^2,\dots,l_n^2\ts|\ts a)$;
see Sec.~\ref{subsec:cps}.

These results imply that the Casimir elements \eqref{genb}
respectively coincide, up to a constant factor, with those found previously in \cite{mn:ci}
as an application of the twisted Yangians. However, the expressions
for those Casimir elements given in \cite{mn:ci} are quite different from
\eqref{genb} so the coincidence appears to be surprising.

Our proofs are based on the characterization theorem for the factorial
Schur polynomials (see \cite{o:qi} and \cite{oo:ss}) as well as on the
eigenvalues of the Jucys--Murphy elements for the Brauer algebra
in the irreducible representations found in \cite{lr:rh} and \cite{n:yo}.
This approach extends to the Casimir elements of the form
\eqref{tracegn}, where $C$ is the anti-symmetrizer in the Brauer algebra.
Although such elements are studied in the literature and their
Harish-Chandra images are known (see e.g. \cite{hu:ci}, \cite{i:td}, \cite{i:tp},
\cite{mn:ci} and \cite{w:ce}), our arguments appear to be new and they apply uniformly
to both families of elements constructed
with the use of the symmetrizers and anti-symmetrizers.

\medskip

We acknowledge the
support of the Australian Research Council.
N.\,I. and E.\,R. are grateful to the University of Sydney
for the warm hospitality during their visits.
N.\,I. was also partially supported by  the Program of Fundamental Research of the Physics
and Astronomy Division of NASU, and Joint Ukrainian-Russian
SFFR-RFBR project F40.2/108.

\section{Brauer algebra}
\label{sec:bra}
\setcounter{equation}{0}

We let $m$ be a positive integer and $\om$ an indeterminate.
An $m$-diagram $d$ is a collection
of $2m$ dots arranged into two rows with $m$ dots in each row
connected by $m$ edges such that any dot belongs to only one edge.
The product of two diagrams $d_1$ and $d_2$ is determined by
placing $d_1$ above $d_2$ and identifying the vertices
of the bottom row of $d_1$ with the corresponding
vertices in the top row of $d_2$. Let $s$ be the number of
closed loops obtained in this placement. The product $d_1d_2$ is given by
$\om^{\tss s}$ times the resulting diagram without loops.
The {\it Brauer algebra\/} $\Bc_m(\om)$ \cite{b:aw} is defined as the
$\CC(\om)$-linear span of the $m$-diagrams with the
multiplication defined above.
The dimension of the algebra is $1\cdot 3\cdots (2m-1)$.
For $1\leqslant a<b\leqslant m$ denote by $s_{a\tss b}$
and $\ep_{a\tss b}$ the respective diagrams
of the form

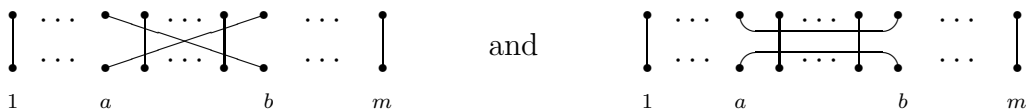
\begin{center}
\begin{picture}(400,60)
\thinlines

\put(10,20){\circle*{3}}
\put(45,20){\circle*{3}}
\put(60,20){\circle*{3}}
\put(90,20){\circle*{3}}
\put(105,20){\circle*{3}}
\put(150,20){\circle*{3}}

\put(10,40){\circle*{3}}
\put(45,40){\circle*{3}}
\put(60,40){\circle*{3}}
\put(90,40){\circle*{3}}
\put(105,40){\circle*{3}}
\put(150,40){\circle*{3}}

\put(10,20){\line(0,1){20}}
\put(60,20){\line(0,1){20}}
\put(45,20){\line(3,1){60}}
\put(45,40){\line(3,-1){60}}
\put(90,20){\line(0,1){20}}
\put(150,20){\line(0,1){20}}

\put(20,20){$\cdots$}
\put(20,35){$\cdots$}
\put(68,20){$\cdots$}
\put(68,35){$\cdots$}
\put(120,20){$\cdots$}
\put(120,35){$\cdots$}

\put(8,5){\scriptsize $1$ }
\put(43,5){\scriptsize $a$ }
\put(105,5){\scriptsize $b$ }
\put(146,5){\scriptsize $m$ }

\put(190,25){\text{and}}

\put(250,20){\circle*{3}}
\put(285,20){\circle*{3}}
\put(300,20){\circle*{3}}
\put(330,20){\circle*{3}}
\put(345,20){\circle*{3}}
\put(390,20){\circle*{3}}

\put(250,40){\circle*{3}}
\put(285,40){\circle*{3}}
\put(300,40){\circle*{3}}
\put(330,40){\circle*{3}}
\put(345,40){\circle*{3}}
\put(390,40){\circle*{3}}

\put(250,20){\line(0,1){20}}
\put(300,20){\line(0,1){20}}
\put(315,20){\oval(60,12)[t]}
\put(315,40){\oval(60,12)[b]}
\put(330,20){\line(0,1){20}}
\put(390,20){\line(0,1){20}}

\put(260,20){$\cdots$}
\put(260,35){$\cdots$}
\put(308,20){$\cdots$}
\put(308,35){$\cdots$}
\put(360,20){$\cdots$}
\put(360,35){$\cdots$}

\put(248,5){\scriptsize $1$ }
\put(283,5){\scriptsize $a$ }
\put(345,5){\scriptsize $b$ }
\put(386,5){\scriptsize $m$ }

\end{picture}
\end{center}

\noindent
The subalgebra of $\Bc_m(\om)$ generated over $\CC$
by $s_{a\ts a+1}$ with $a=1,\dots,m-1$
is isomorphic to the group algebra of the symmetric group
$\CC[\Sym_m]$ so that $s_{a\tss b}$
will be identified with the transposition $(a\ts b)$.
The Brauer algebra $\Bc_{m-1}(\om)$ will be regarded
as a natural subalgebra of $\Bc_m(\om)$.

The {\it Jucys--Murphy elements\/} $x_1,\dots,x_m$
for the Brauer algebra $\Bc_m(\om)$ are given by the formulas
\beql{jmdef}
x_1=0,\qquad x_b=\sum_{a=1}^{b-1}(s_{a\tss b}-\ep_{a\tss b}),
\qquad b=2,\dots,m;
\eeq
see \cite{lr:rh} and \cite{n:yo}, where, in particular,
the eigenvalues
of the $x_b$ in irreducible representations were calculated.
The element $x_m$ commutes
with the subalgebra $\Bc_{m-1}(\om)$. This implies that
the elements $x_1,\dots,x_m$
of $\Bc_m(\om)$ pairwise commute.

Irreducible representations of the algebra $\Bc_m(\om)$ (over $\CC(\om)$)
are parameterized
by the set of partitions of the numbers $m-2f$
with $f\in\{0,1,\dots,\lfloor m/2\rfloor\}$.
We will identify partitions with their diagrams
so that if the parts of $\la$ are $\la_1,\la_2,\dots$ then
the corresponding diagram is
a left-justified array of rows of unit boxes containing
$\la_1$ boxes in the top row, $\la_2$ boxes
in the second row, etc. We will denote by $|\la|$
the number of boxes in the diagram and by
$\ell(\la)$ its {\it length\/}, i.e., the number of rows.
The box in row $i$ and column $j$ of a diagram
will be denoted as the pair $(i,j)$.
An {\it updown $\la$-tableau\/} is a sequence
$T=(\La_1,\dots,\La_m)$ of $m$ diagrams such that
for each $r=1,\dots,m$ the diagram $\La_r$ is obtained
from $\La_{r-1}$ by adding or removing one box,
where $\La_0=\varnothing$ is the empty diagram and $\La_m=\la$.
To each updown tableau $T$ we associate the corresponding
sequence of {\it contents\/} $(c_1,\dots,c_m)$, $c_r=c_r(T)$,
where
\ben
c_r=j-i\qquad\text{or}\qquad
c_r=-\om+1-j+i,
\een
if $\La_r$ is obtained by adding the box $(i,j)$ to $\La_{r-1}$
or by removing this box from $\La_{r-1}$,
respectively.

It is well-known that the Jucys--Murphy elements can be used to define the
primitive idempotents $E_{T}=E^{\lambda}_{T}$; see e.g. \cite{im:fp}
for explicit formulas.
When $\la$ runs over all partitions
of $m,m-2,\dots$ and $T$ runs over all updown $\la$-tableaux,
the elements $\{E_{T}\}$ yield a complete set
of pairwise orthogonal primitive
idempotents for $\Bc_m(\om)$. They have the properties
\beql{xiet}
x_r\ts E_{T}=E_{T}\ts x_r=c_r(T)\ts E_{T}, \qquad r=1,\dots,m;
\eeq
see \cite{lr:rh} and \cite{n:yo}.
In particular, if $\la=(m)$ is the single row diagram with $m$ boxes
then there is a unique updown $\la$-tableau which can also be regarded
as the standard tableau obtained by writing the numbers $1,\dots,m$ into
the boxes of $\la$ from left to right. The corresponding primitive idempotent
is the {\it symmetrizer\/} $S^{(m)}\in\Bc_m(\om)$ given in terms of the
Jucys--Murphy elements as
\beql{symjm}
S^{(m)}=\prod_{r=2}^m\frac{(1+x_r)(\om+r-3+x_r)}{r\tss(\om+2r-4)}.
\eeq
This element also admits an equivalent expression
\beql{symfpnew}
S^{(m)}=\frac{1}{m!}
\prod_{1\leqslant a<b\leqslant m}
\Big(1+\frac{s_{ab}}{b-a}-\frac{\ep_{ab}}
{\om/2+b-a-1}\Big),
\eeq
where the product is taken in the lexicographic order
on the pairs $(a\ts b)$; see \cite{im:fp} and \cite{m:ff} for some
other equivalent formulas for $S^{(m)}$. We have the properties
\beql{mulse}
s_{ab}\tss S^{(m)}=S^{(m)}\tss s_{ab}=S^{(m)}\Fand
\ep_{ab}\tss S^{(m)}=S^{(m)}\tss\ep_{ab}=0
\eeq
for all $1\leqslant a<b\leqslant m$.

Similarly, if $\la=(1^m)$ is the single column diagram with $m$ boxes
then the unique updown $\la$-tableau can be regarded
as the standard tableau obtained by writing the numbers $1,\dots,m$ into
the boxes of $\la$ from top to bottom. The corresponding primitive idempotent
is the {\it anti-symmetrizer\/} $A^{(m)}\in\Bc_m(\om)$.
It is well-known that $A^{(m)}$ coincides with the anti-symmetrizer
in the group algebra for the symmetric group $\Sym_m$ and has
the properties
\beql{mulant}
s_{ab}\tss A^{(m)}=A^{(m)}\tss s_{ab}=-A^{(m)}\Fand
\ep_{ab}\tss A^{(m)}=A^{(m)}\tss\ep_{ab}=0
\eeq
for all $1\leqslant a<b\leqslant m$.

Consider now the action of the Brauer algebra
on the tensor space \eqref{tenprk}. In the orthogonal case, $\om=N$
and the generators of $\Bc_m(N)$ act by the rule
\beql{braact}
s_{ab}\mapsto P_{ab},\qquad \ep_{ab}\mapsto Q_{ab},\qquad
a< b,
\eeq
where $P_{ab}$ is defined by
\beql{pdef}
P_{ab}=\sum_{i,j=1}^N 1^{\ot(a-1)}\ot e_{ij}
\ot 1^{\ot(b-a-1)}\ot e_{ji}\ot 1^{\ot(m-b)},
\eeq
while
\beql{qdefo}
Q_{ab}=\sum_{i,j=1}^N 1^{\ot(a-1)}\ot e_{ij}
\ot 1^{\ot(b-a-1)}\ot e_{i'j'}\ot 1^{\ot(m-b)}
\eeq
and we use the involution on the set
of indices $\{1,\dots,N\}$ defined by $i\mapsto i^{\tss\prime}=N-i+1$.

In the symplectic case, $\om=-N$ (where $N=2n$ is even) and the action of
$\Bc_m(-N)$ in the space \eqref{tenprk}
is now defined by
\beql{braactsp}
s_{ab}\mapsto -P_{ab},\qquad \ep_{ab}\mapsto -Q_{ab},\qquad
a< b,
\eeq
where $P_{ab}$ is defined in \eqref{pdef}, and
\beql{qdefsp}
Q_{ab}=\sum_{i,j=1}^N \ve_i\tss\ve_j\ts 1^{\ot(a-1)}\ot e_{ij}
\ot 1^{\ot(b-a-1)}\ot e_{i'j'}\ot 1^{\ot(m-b)}
\eeq
with $\ve_i=1$ for $i=1,\dots,n$ and
$\ve_i=-1$ for $i=n+1,\dots,2n$.

\section{Harish-Chandra images}
\label{sec:cas}
\setcounter{equation}{0}

Consider the Lie algebra $\gl_N$ with its standard basis
elements $E_{ij}$, $1\leqslant i,j\leqslant N$.
The Lie subalgebra
of $\gl_N$ spanned by the elements $F_{ij}=E_{ij}-E_{j'i'}$ is isomorphic to
the orthogonal Lie algebra $\oa_N$. Similarly, the Lie subalgebra
of $\gl_{2n}$ spanned by the elements
$F_{ij}=E_{ij}-\ve_i\tss\ve_j\tss E_{j'i'}$ is isomorphic to
the symplectic Lie algebra $\spa_{2n}$.
We will keep the notation $\g_N$ for the Lie algebra $\oa_N$
(with $N=2n$ or $N=2n+1$) or $\spa_N$ (with $N=2n$).

Using the notation \eqref{ea} we can write
the defining relations of the universal enveloping algebra $\U(\g_N)$
in the matrix form as
\beql{deff}
F_1\ts F_2-F_2\ts F_1=(P-Q)\ts F_2-F_2\ts (P-Q)
\eeq
together with the relation $F+F^{\tss\prime}=0$, where the prime denotes
the matrix transposition defined by
\beql{transp}
(A^{\tss\prime})_{ij}=\begin{cases} A_{j'i'}\qquad&\text{in the orthogonal case,}\\
\ve_i\tss\ve_j\ts A_{j'i'}\qquad&\text{in the symplectic case.}
\end{cases}
\eeq

In the following lemma we identify the elements of the respective
Brauer algebra $\Bc_m(N)$ or $\Bc_m(-N)$ with their images
under the actions \eqref{braact} or \eqref{braactsp}.

\ble\label{lem:invpe}
Let $u_1,\dots,u_m$ be complex parameters.
For any permutations $\si,\tau\in\Sym_m$ and for $C=S^{(m)}$ or $C=A^{(m)}$ we have
\beql{tracegnp}
\tr\ts (F_{\si(1)}+u_{\tau(1)})\dots (F_{\si(m)}+u_{\tau(m)})\ts C=
\tr\ts (F_1+u_1)\dots (F_m+u_m)\ts C.
\eeq
\ele

\bpf
Let $P_{\pi}$ denote the image of any element $\pi\in\Sym_m$ in the
algebra \eqref{tenprku}. Since $P_{\pi}\tss F_a=F_{\pi(a)}\tss P_{\pi}$,
using the cyclic property of trace
and applying the conjugation by the element $P_{\si^{-1}}$ in the
left hand side of \eqref{tracegnp} we find that it suffices to
verify the relation for the case where $\si$ is the identity permutation.
Note that by \eqref{deff},
\begin{multline}
(F_a+u)(F_{a+1}+v)-(F_{a+1}+v)(F_a+u)\\
=(P_{a\ts a+1}-Q_{a\ts a+1})\tss F_{a+1}
-F_{a+1}\tss(P_{a\ts a+1}-Q_{a\ts a+1}),
\non
\end{multline}
and by \eqref{mulse} and \eqref{mulant}
\ben
C\ts(P_{a\ts a+1}-Q_{a\ts a+1})=(P_{a\ts a+1}-Q_{a\ts a+1})\ts C=\pm\tss C.
\een
Hence, the claim follows from the cyclic property of trace and the
first part of the proof.
\epf

Given any $n$-tuple of complex numbers
$\la=(\la_1,\dots,\la_n)$ the corresponding irreducible highest weight
representation $L(\la)$ of the Lie algebra $\g_N$ is generated by
a nonzero vector $\xi\in L(\la)$ such that
\begin{alignat}{2}
F_{ij}\ts\xi&=0 \qquad &&\text{for} \quad
1\leqslant i<j\leqslant N, \qquad \text{and}
\non
\\
F_{ii}\ts\xi&=\la_i\ts\xi \qquad &&\text{for} \quad 1\leqslant i\leqslant n.
\non
\end{alignat}

We will denote by $\Gr_N$ the orthogonal group $\Or_N$ or the symplectic group $\Spr_N$.
Recall that finite-dimensional irreducible representations of the orthogonal group
$\Or_N$ are parameterized by all diagrams $\la$ with the property
$\la'_1+\la'_2\leqslant N$, where $\la'_j$ denotes the number
of boxes in the column $j$ of $\la$.
The corresponding representation will be denoted by $V(\la)$.
Let $\la^*$ be
the diagram obtained from $\la$ by replacing the first column
with the column containing $N-\la'_1$ boxes.
If $N=2n+1$ and $\la'_1\leqslant n$ then
the associated representation
of the Lie algebra $\oa_N$ in the space $V(\la)$
is irreducible and isomorphic to the representation $L(\la)$ whose highest weight
coincides with $\la$; if $\la'_1> n$ then
the associated representation
of $\oa_N$ is isomorphic to $L(\la^*)$.

If $N=2n$ and $\la'_1< n$ then
the associated representation
of the Lie algebra $\oa_N$ in the space $V(\la)$
is irreducible and isomorphic to $L(\la)$, while for $\la'_1> n$ the
associated representation
of $\oa_N$ is isomorphic to $L(\la^*)$. If $N=2n$ and $\la'_1=n$
then the associated representation
of $\oa_N$ in $V(\la)$ is isomorphic to the direct sum
of two irreducible representations $L(\la)$ and $L(\wt\la)$
with $\wt\la=(\la_1,\dots,\la_{n-1},-\la_n)$.

Finite-dimensional irreducible representations $V(\la)$ of the symplectic group
$\Spr_N$ with $N=2n$ are parameterized by partitions $\la$
whose lengths do not exceed $n$. The associated representation of the Lie algebra
$\spa_N$ in $V(\la)$ is irreducible and isomorphic to $L(\la)$.

Any element $z\in\U(\g_N)^{\Gr_N}$ acts in $V(\la)$
by multiplying each vector by a scalar $\chi(z)$.
In the case of $\g_N=\oa_{2n}$ the eigenvalues of any element
$z\in\U(\oa_{2n})^{\Or_{2n}}$ in $L(\la)$ and $L(\wt\la)$ coincide and they are equal to
$\chi(z)$.

When regarded as a function of the highest weight, $\chi(z)$
is a symmetric polynomial in the variables $l_1^2,\dots,l_n^2$,
where $l_i=\la_i+\rho_i$ and $\rho_i=n-i+\ve$ with
\beql{rhocc}
\ve=\begin{cases}
0\quad&\text{for}\quad\g_N=\oa^{}_{2n},\\
\frac12\quad&\text{for}\quad\g_N=\oa^{}_{2n+1},\\
1\quad&\text{for}\quad\g_N=\spa^{}_{2n}.
\end{cases}
\eeq
The mapping $z\mapsto\chi(z)$ defines an algebra isomorphism
\ben
\chi:\U(\g_N)^{\Gr_N}\to\CC[l_1^2,\dots,l_n^2]^{\Sym_n}
\een
known as the Harish-Chandra isomorphism; see e.g. \cite[Ch.~7]{d:ae}.

\subsection{Characterization properties for symmetric polynomials}
\label{subsec:cps}

Following \cite{mn:ci},
we will be using the factorial (or double) elementary
and complete symmetric polynomials and their characterization properties;
see \cite{o:qi} and \cite{oo:ss} for more details. Here we recall the corresponding
results to be used in the proofs below.

Consider the algebra of symmetric polynomials in the
independent variables $z_1,\dots, z_n$ over $\CC$ and fix a sequence
$a=(a_1,a_2,\dots)$ of complex numbers.
The {\it factorial elementary {\rm and} complete
symmetric polynomials\/} are defined by the respective formulas
\begin{align}\label{felem}
e_k(z_1,\dots, z_n\tss|\tss a)&
=\sum_{1\leqslant p_1<\dots <p_k\leqslant n}(z_{p_1}-a_{p_1})
(z_{p_2}-a_{p_2-1})\dots (z_{p_k}-a_{p_k-k+1}),\\
\label{fcoml}
h_k(z_1,\dots, z_n\tss|\tss a)&
=\sum_{1\leqslant p_1\leqslant\dots\leqslant p_k\leqslant n}(z_{p_1}-a_{p_1})
(z_{p_2}-a_{p_2+1})\dots (z_{p_k}-a_{p_k+k-1}),
\end{align}
so that $e_k(z_1,\dots, z_n\tss|\tss a)=0$ for $k>n$.
These polynomials are particular cases of the factorial (or double) Schur polynomials;
see e.g. \cite{m:sft}. When $a$ is specialized to the sequence
of zeros, then \eqref{felem} and \eqref{fcoml} become the elementary and complete
symmetric polynomials $e_k(z_1,\dots, z_n)$ and $h_k(z_1,\dots, z_n)$.

For any partition $\la=(\la_1,\dots,\la_n)$ whose length $\ell(\la)$ does not
exceed $n$ introduce
the $n$-tuple $a_{\la}$ of complex numbers by
\ben
a_{\la}=(a_{\la_1+n},a_{\la_2+n-1},\dots,a_{\la_n+1}).
\een
The polynomials \eqref{felem}
and \eqref{fcoml} possess vanishing properties of the form:
\begin{alignat}{2}
&\text{if}\quad\ell(\la)<k\qquad &&\text{then}\quad e_k(a_{\la}\tss|\tss a)=0,
\non\\
&\text{if}\quad\la_1<k\qquad &&\text{then}\quad h_k(a_{\la}\tss|\tss a)=0.
\non
\end{alignat}
We will need two particular cases of the characterization theorem for
the factorial Schur polynomials \cite{o:qi}. Now we will be assuming that
all elements $a_i$ of the sequence $a$ are distinct.
Suppose that $f(z_1,\dots,z_n)$ is
a symmetric polynomial of degree $\leqslant k$
whose component of degree $k$ coincides with
$e_k(z_1,\dots, z_n)$ or $h_k(z_1,\dots, z_n)$.
If
\ben
f(a_{\la})=0\qquad\text{for all $\la$ with $|\la|<k$}
\een
then $f(z_1,\dots,z_n)$ equals
$e_k(z_1,\dots, z_n\tss|\tss a)$
or $h_k(z_1,\dots, z_n\tss|\tss a)$, respectively.

From now on, we will work with particular sequences $a$ defined by
\beql{seqa}
a=(\ve^2,(\ve+1)^2,(\ve+2)^2,\dots),
\eeq
where $\ve$ is introduced in \eqref{rhocc}, so that $a_i=(\ve+i-1)^2$.
Furthermore, the $n$-tuple $a_{\la}$ associated with the sequence \eqref{seqa}
has the form
\ben
a_{\la}=\big((\la_1+n-1+\ve)^2,\dots,(\la_n-1+\ve)^2\big)=(l_1^2,\dots,l_n^2).
\een
Note that any element $z\in\U(\g_N)^{\Gr_N}$ is uniquely determined by
the eigenvalues $\chi(z)$ in the irreducible modules $L(\la)$, where
$\la=(\la_1,\dots,\la_n)$ runs over the set of partitions with $\ell(\la)\leqslant n$.
Hence, we come to the following characterization properties
of Casimir elements; cf. \cite[Corollary~2.5]{mn:ci}.
We use the canonical filtration on the universal enveloping algebra $\U(\g_N)$.

\bpr\label{prop:charz}
Suppose that $z\in\U(\g_N)^{\Gr_N}$ is an element of degree $\leqslant 2k$
which vanishes in each representation $L(\la)$
with $|\la|<k$. If the homogeneous component of $\chi(z)$
of degree $2k$ coincides with $e_k(\la^2_1,\dots, \la^2_n)$
or $h_k(\la^2_1,\dots, \la^2_n)$ then $\chi(z)$ equals
$e_k(l^2_1,\dots, l^2_n\tss|\tss a)$
or $h_k(l^2_1,\dots, l^2_n\tss|\tss a)$, respectively.
\epr

\subsection{Casimir elements for the orthogonal Lie algebras}
\label{subsec:ceo}

As we pointed out in the Introduction, any element of the universal
enveloping algebra $\U(\g_N)$ of the form \eqref{tracegn} belongs to
the subalgebra of invariants $\U(\g_N)^{\Gr_N}$. We will be concerned
with two choices of the element $C\in\Bc_m(\om)$; namely,
$C=S^{(m)}$ and $C=A^{(m)}$. Moreover, the parameters $u_a$ in
\eqref{tracegn} will be specialized accordingly so that all differences
$u_a-u_{a+1}$ have the same value $1$ or $-1$ for all $a=1,\dots,m-1$.

We will be assuming here that $\g_N=\oa_N$. For any $m\geqslant 0$ set
\beql{alpm}
\al_m=\frac{N+2m-2}{N+m-2}.
\eeq
We start by taking even values $m=2k$ and particular specializations
of the parameters.

\bth\label{thm:caseven}
For any $k\geqslant 1$ the image of the Casimir element
\beql{genbth}
\tr\ts (F_1+k-1)\dots (F_{2k}-k)\ts S^{(2k)}\in \U(\oa_N)^{\Or_N}
\eeq
under the Harish-Chandra isomorphism coincides with
$\al_{2k}\ts h_k(l^2_1,\dots, l^2_n\tss|\tss a)$.
\eth

\bpf
Denote the Casimir element \eqref{genbth} by $D_k$. We will use
Proposition~\ref{prop:charz} and start by showing that $D_k$ vanishes
in all representations $L(\la)$ of $\oa_N$, where the partitions
$\la$ satisfy $|\la|<k$. We employ a realization of $L(\la)$ in tensor
spaces as follows. Let $r=|\la|$. Consider the action of the Lie algebra $\oa_N$
on $\CC^N$ defined by
\ben
F_{ij}\mapsto -e_{ji}+e_{i'j'},\qquad 1\leqslant i,j\leqslant N,
\een
so that this representation is isomorphic to $L(1,0,\dots,0)$.
The space of tensors $(\CC^N)^{\ot r}$ then also becomes a representation of $\oa_N$.
Using the matrix notation \eqref{ea}, under the corresponding homomorphism
\ben
\vp:\U(\oa_N)\ot\End(\CC^N)^{\ot 2k}\to \End(\CC^N)^{\ot r}\ot \End(\CC^N)^{\ot 2k}
\een
we have
\beql{imfa}
\vp(F_a)= \sum_{b=1}^r (-P_{\tss b\ts\ts r+a}+Q_{\tss b\ts\ts r+a}),\qquad a=1,\dots,2k,
\eeq
where we use the operators \eqref{pdef} and \eqref{qdefo}.
The decomposition of $(\CC^N)^{\ot r}$ into a direct sum of irreducible
representations of $\oa_N$ contains $L(\la)$ with a nonzero multiplicity.
Hence, the desired vanishing condition of Proposition~\ref{prop:charz}
will follow if we show that
\beql{etfze}
\big(\vp(F_1)+k-1\big)\dots \big(\vp(F_{2k})-k\big)\ts S^{\tss\prime(2k)}=0,
\eeq
where $S^{\tss\prime(2k)}$ denotes the image of the symmetrizer in the Brauer
algebra $\Bc_{2k}(N)$ acting of the last $2k$ copies of the tensor product space
$(\CC^N)^{\ot(r+2k)}$. Due to \eqref{imfa} and relations \eqref{mulse},
the desired identity \eqref{etfze} can be written in the form
\beql{etjm}
(-x_{r+1}+k-1)\dots (-x_{r+2k}+k-1)\ts S^{\tss\prime(2k)}=0,
\eeq
where $x_{r+1},\dots,x_{r+2k}$ denote the images of the Jucys--Murphy
elements \eqref{jmdef} of the Brauer algebra $\Bc_{r+2k}(N)$
under its action in $(\CC^N)^{\ot(r+2k)}$.
To prove \eqref{etjm}
we note that the two operators
$(-x_{r+1}+k-1)\dots (-x_{r+2k}+k-1)$
and $S^{\tss\prime(2k)}$ on the vector space $(\CC^N)^{\ot(r+2k)}$
commute with the action of $\Or_N$
and show that their images
have zero intersection. To describe the image of the first operator,
represent the vector space as the direct sum of irreducible
representations of $\Or_N$,
\ben
(\CC^N)^{\ot(r+2k)}=\sum_{l=0}^{\lfloor \frac{r}{2}\rfloor+k}
\sum_{\nu\vdash r+2k-2l}\sum_U E_U(\CC^N)^{\ot(r+2k)},
\een
where the last sum is taken over all updown tableaux
$U=(\La_{1},\dots,\La_{r+2k})$ of shape $\nu$
associated with the Brauer algebra $\Bc_{r+2k}(N)$.
We claim that if $U$ is an updown
tableau of shape $\nu=\La_{r+2k}$ with $\nu_1\geqslant k$ then
\beql{xeu}
(-x_{r+1}+k-1)\dots (-x_{r+2k}+k-1)\ts E_U=0.
\eeq
Indeed, since $r<k$ there exists
a pair
of diagrams $(\La_{r+a},\La_{r+a+1})$ with $a\in\{0,1,\dots,2k-1\}$
with the property that the second diagram is obtained from the first
by adding the box $(1,k)$. The content of this box is $k-1$ so that
\eqref{xeu} follows from relations \eqref{xiet}. Thus, as a representation
of $\Or_N$ the image
\ben
(-x_{r+1}+k-1)\dots (-x_{r+2k}+k-1)\ts (\CC^N)^{\ot(r+2k)}
\een
is contained in a direct sum of representations $V(\nu)$ with $\nu_1<k$.

On the other hand, the operator $S^{\tss\prime(2k)}$
projects the vector space $(\CC^N)^{\ot(r+2k)}$ onto the tensor
product $(\CC^N)^{\ot r}\ot V(2k,0,\dots,0)$ of representations of $\Or_N$.
For the irreducible decomposition of the tensor product
of representations of $\Or_N$ we have
\ben
\CC^N\ot V(\mu)\cong \bigoplus_{\wt\mu}\ts V(\tss\wt\mu\tss),
\een
where $\wt\mu$ is obtained from $\mu$ by adding or removing one box.
Hence, as a representation
of $\Or_N$, the image
$S^{\tss\prime(2k)}(\CC^N)^{\ot(r+2k)}$
is contained in a direct sum of representations $V(\nu)$ with $\nu_1\geqslant k+1$.
This completes the proof of \eqref{etjm} and hence \eqref{etfze}.

The leading term of the symmetric
polynomial $\chi(D_k)$ was calculated in \cite{m:ff}
and \cite{mr:cm}. It coincides with $\al_{2k}\tss h_k(\la^2_1,\dots, \la^2_n)$
so that the proof is completed by the application of
Proposition~\ref{prop:charz}.
\epf

\bre\label{rem:twist} (i)\quad
Theorem~\ref{thm:caseven} together with Theorem~\ref{thm:scaseven} below show that the elements
\eqref{genbth} and \eqref{sgenbth} are proportional to the respective Casimir elements
given in \cite[Theorems~3.2 and~3.3]{mn:ci} by completely different formulas.
It would be interesting to find a direct argument connecting these families.
\par
(ii)\quad It is possible to give an independent proof of
Theorems~\ref{thm:caseven} without relying
on the calculation of the leading term of the symmetric
polynomial $\chi(D_k)$ in \cite{m:ff}
and \cite{mr:cm}. To this end, we could use a different
version of the characterization theorem from \cite{o:qi}
where the assumption $|\la|<k$ is replaced by $\la_1<k$.
The vanishing condition is verified in a similar way and it determines
the polynomial $\chi(D_k)$, up to a constant factor. The latter
is calculated by finding the leading term in the case where
$\la$ is the one box diagram.
\qed
\ere

Now we turn to more general Casimir elements of the form \eqref{tracegn}
and \eqref{tracedi}.
Let $u$ be a variable. Consider the polynomials
in $u$ whose coefficients are Casimir elements for $\oa_N$ given by
\ben
D_{m}(u)=
\tr\tss \Big(F_1+u+\frac{m-1}{2}\Big)\Big(F_2+u+\frac{m-3}{2}\Big)
\dots \Big(F_{m}+u-\frac{m-1}{2}\Big)\ts S^{(m)}.
\een
We use notation \eqref{alpm}.

\bco\label{cor:dpol}
For the images
under the Harish-Chandra isomorphism we have
\ben
\chi:D_{m}(u)\mapsto \al_{m}
\ts\sum_{r=0}^{\lfloor\frac{m}{2}\rfloor}\ts
\binom{N+m-2}{m-2r}\ts h_{r}(l^2_1,\dots, l^2_n\tss|\tss a)
\ts \prod_{i=0}^{m-2r-1}\Big(u-\frac{m-1}{2}+r+i\Big)
\een
and
\begin{multline}
\chi:\tr\ts (-\di_t+F_1\tss t^{-1})\dots (-\di_t+F_m\tss t^{-1})\ts S^{(m)}
\mapsto \\
\al_{m}\ts\sum_{r=0}^{\lfloor\frac{m}{2}\rfloor}\ts
\binom{N+m-2}{m-2r}\ts h_{r}(l^2_1,\dots, l^2_n\tss|\tss a)
\ts t^{-2r}\ts (-\di_t+r\tss t^{-1})^{m-2r}.
\non
\end{multline}
\eco

\bpf
Observe that by the property of trace,
$D_m(u)$ is stable under
the transposition \eqref{transp} applied simultaneously
to each of the $m$ copies
of the algebra $\End\CC^N$. The symmetrizer $S^{(m)}$ is also stable
under this transposition. On the other hand,
since $F^{\tss\prime}=-F$, using Lemma~\ref{lem:invpe},
for the image of $D_m(u)$ we find
\ben
\bal
D_{m}(u)&=
\tr\tss {\Big(-}F_1+u+\frac{m-1}{2}\Big)
\dots {\Big(-}F_{m}+u-\frac{m-1}{2}\Big)\ts S^{(m)}\\
{}&=(-1)^m\ts\tr\tss \Big(F_1-u-\frac{m-1}{2}\Big)
\dots \Big(F_{m}-u+\frac{m-1}{2}\Big)\ts S^{(m)}\\
{}&=(-1)^m\ts\tr\tss \Big(F_1-u+\frac{m-1}{2}\Big)
\dots \Big(F_{m}-u-\frac{m-1}{2}\Big)\ts S^{(m)}=(-1)^m\ts D_{m}(-u).
\eal
\een
This shows that the polynomials $D_{2k}(u)$ are even, while
the polynomials $D_{2k-1}(u)$ are odd.
In particular, $D_{2k-1}(0)=0$, while the value
$\chi(D_{2k}(1/2))=\chi(D_{2k}(-1/2))$ is found
from Theorem~\ref{thm:caseven}. These values agree with the
Harish-Chandra images provided in the statement of the corollary.
Hence, the proof will be completed if we show that the polynomials $D_m(u)$
and the polynomials which are claimed to be their Harish-Chandra images
satisfy the same recurrence relations. We show first that
\beql{drecr}
D_m(u+1/2)-D_m(u-1/2)=\frac{(N+m-3)(N+2m-2)}{N+2m-4}\ts D_{m-1}(u),\qquad m\geqslant 1,
\eeq
where we set $D_0(u)=1$. By Lemma~\ref{lem:invpe},
\ben
D_m(u+1/2)=\tr\tss \Big(F_1+u+\frac{m}{2}-1\Big)
\dots \Big(F_{m-1}+u-\frac{m}{2}-1\Big)\ts
\Big(F_{m}+u+\frac{m}{2}\Big)\ts S^{(m)},
\een
so that
\ben
D_m(u+1/2)-D_m(u-1/2)=m\ts \tr
\ts\Big(F_1+u+\frac{m}{2}-1\Big)
\dots \Big(F_{m-1}+u-\frac{m}{2}-1\Big)\ts S^{(m)}.
\een
By \cite[Lemma~4.1]{m:ff}, the partial trace of the symmetrizer is
found by
\ben
\tr_m\ts S^{(m)}=\frac{(N+m-3)(N+2m-2)}{m\ts(N+2m-4)}\ts
S^{(m-1)}
\een
thus verifying \eqref{drecr}. A simple calculation shows that
the same relation is satisfied by
the polynomials which are claimed to be the images $\chi(D_m(u))$ as stated
in the corollary.

To prove the second part, note that by the relation
\ben
t^m\ts (-\di_t+F_1\tss t^{-1})\dots (-\di_t+F_m\tss t^{-1})
=(-t\tss\di_t+F_1+m-1)\dots (-t\tss\di_t+F_m)
\een
the polynomial $D_m(u)$ can be written in the form
\ben
t^m\ts\tr\ts (-\di_t+F_1\tss t^{-1})\dots (-\di_t+F_m\tss t^{-1})\ts S^{(m)}
\een
after the subsequent replacement of $-t\tss\di_t$ with $u-(m-1)/2$.
Similarly, for any $k\geqslant 0$,
\ben
t^k\ts (-\di_t+r\tss t^{-1})^k=(-t\tss\di_t+r+k-1)\dots (-t\tss\di_t+r),
\een
so that the second relation follows from the first.
\epf

In what follows we state analogues of Theorem~\ref{thm:caseven}
and Corollary~\ref{cor:dpol}, where
the role of the symmetrizer $S^{(m)}$ is taken by
the anti-symmetrizer $A^{(m)}$. The arguments are quite similar
so we only indicate the key steps in the proofs. The corresponding
Casimir elements turn out to coincide with those already appeared
in the literature; cf. e.g. \cite{i:td}, \cite{mn:ci} and \cite{w:ce}.

\bth\label{thm:ascaseven}
For any $1\leqslant k\leqslant n$ the image of the Casimir element
\beql{asgenbth}
\tr\ts (F_1-k+1)\dots (F_{2k}+k)\ts A^{(2k)}\in \U(\oa_N)^{\Or_N}
\eeq
under the Harish-Chandra isomorphism coincides with
$(-1)^k\ts e_k(l^2_1,\dots, l^2_n\tss|\tss a)$.
\eth

\bpf
Denote the Casimir element \eqref{asgenbth} by $C_k$. We use
Proposition~\ref{prop:charz} and show that $C_k$ vanishes
in all representations $L(\la)$ of $\oa_N$, where the partitions
$\la$ satisfy $|\la|<k$. Using the same realization of $L(\la)$
as in the proof of Theorem~\ref{thm:caseven}, we come to showing
the following analogue of \eqref{etfze}:
\beql{asetfze}
(\vp(F_1)-k+1\big)\dots \big(\vp(F_{2k})+k\big)\ts A^{\tss\prime(2k)}=0.
\eeq
Here
the homomorphism $\vp$ is defined in \eqref{imfa} and
$A^{\tss\prime(2k)}$ denotes the image of the anti-symmetrizer in the Brauer
algebra $\Bc_{2k}(N)$ acting of the last $2k$ copies of the tensor product space
$(\CC^N)^{\ot(r+2k)}$ with $r=|\la|$.
By \eqref{mulant},
the identity \eqref{asetfze} can be written in the form
\ben
(-x_{r+1}-k+1)\dots (-x_{r+2k}-k+1)\ts A^{\tss\prime(2k)}=0.
\een
To verify this relation,
we show exactly as in the proof of Theorem~\ref{thm:caseven} that
the image of the operator
$(-x_{r+1}-k+1)\dots (-x_{r+2k}-k+1)$
on the vector space $(\CC^N)^{\ot(r+2k)}$
is contained in a direct sum of representations $V(\nu)$ of $\Or_N$
with $\ell(\nu)<k$,
while the image of the operator $A^{\tss\prime(2k)}$
is contained in a direct sum of representations $V(\nu)$ with $\ell(\nu)\geqslant k+1$.

It was shown in \cite{m:ff}
and \cite{mr:cm} that
the leading term of the symmetric
polynomial $\chi(C_k)$ coincides with $(-1)^k\tss e_k(\la^2_1,\dots, \la^2_n)$.
\epf

Consider the polynomials
in $u$ whose coefficients are Casimir elements for $\oa_N$ given by
\ben
C_{m}(u)=
\tr\tss \Big(F_1+u+\frac{m-1}{2}\Big)\Big(F_2+u+\frac{m-3}{2}\Big)
\dots \Big(F_{m}+u-\frac{m-1}{2}\Big)\ts A^{(m)}.
\een

\bco\label{cor:asdpol}
For the images
under the Harish-Chandra isomorphism we have
\ben
\chi:C_{m}(u)\mapsto
\sum_{r=0}^{\lfloor\frac{m}{2}\rfloor}\ts(-1)^r\ts
\binom{N-2r}{m-2r}\ts
e_{r}(l^2_1,\dots, l^2_n\tss|\tss a)
\ts \prod_{i=0}^{m-2r-1}\Big(u-\frac{m-1}{2}+r+i\Big)
\een
and
\begin{multline}
\chi:\tr\ts (-\di_t+F_1\tss t^{-1})
\dots (-\di_t+F_m\tss t^{-1})\ts A^{(m)}\mapsto\\
\sum_{r=0}^{\lfloor\frac{m}{2}\rfloor}\ts(-1)^r\ts
\binom{N-2r}{m-2r}\ts
e_{r}(l^2_1,\dots, l^2_n\tss|\tss a)
\ts t^{-2r}\ts (-\di_t+r\tss t^{-1})^{m-2r}.
\non
\end{multline}
\eco

\bpf
We argue as in the proof of Corollary~\ref{cor:dpol}.
Using Lemma~\ref{lem:invpe}
we show first that $C_{m}(u)=(-1)^m\ts C_{m}(-u)$.
In particular, $C_{2k-1}(0)=0$ and
$\chi(C_{2k}(1/2))=\chi(C_{2k}(-1/2))$ is found
from Theorem~\ref{thm:ascaseven}. As the next step,
we verify that
\beql{asdrecr}
C_m(u+1/2)-C_m(u-1/2)=(N-m+1)\ts C_{m-1}(u),\qquad m\geqslant 1,
\eeq
where $C_0(u)=1$. This follows easily from Lemma~\ref{lem:invpe},
and the calculation
of the partial trace of the anti-symmetrizer which is
found by
\ben
\tr_m\ts A^{(m)}=\frac{N-m+1}{m}\ts
A^{(m-1)}
\een
thus verifying \eqref{asdrecr}. The same relation is satisfied by
the polynomials which are claimed to be the images $\chi(C_m(u))$ as stated
in the corollary. The second part follows from the first; see
the proof of Corollary~\ref{cor:dpol}.
\epf

\subsection{Casimir elements for the symplectic Lie algebras}
\label{subsec:sceo}

Now take $\g_N=\spa_{2n}$.
Consider the action of the Brauer algebra $\Bc_m(-2n)$
on the space \eqref{tenprk} defined by \eqref{braactsp}.
Note that the image of the symmetrizer $S^{(m)}$ (see
\eqref{symjm} and \eqref{symfpnew}) under this action
is well-defined for $m\leqslant n+1$ and it is zero
for $m=n+1$, while the specialization of
$S^{(m)}$ at $\om=-2n$ is not defined for
$n+2\leqslant m\leqslant 2n$. Nevertheless,
the expression
\beql{deftrsp}
\frac{1}{n-m+1}\ts
\tr\ts S^{(m)} (F_1+u_1)\dots (F_m+u_m)
\eeq
still defines a Casimir element for $\spa_{2n}$ for all
$1\leqslant m\leqslant 2n$, where $u_1,\dots,u_m$ are arbitrary
complex numbers. Indeed, assuming that $m\leqslant n$,
using \eqref{symjm} and calculating
first the partial trace $\tr_m$ in \eqref{deftrsp} over the $m$-th copy
of $\End\CC^{2n}$ we get an expression involving the symmetrizer $S^{(m-1)}$
with the extra factor $(n-m+1)/(n-m+2)$. The latter expression is well-defined
for $m\leqslant n+1$ thus allowing us to extend the value of \eqref{deftrsp}
to $m=n+1$. Continuing with a similar calculation and taking further partial traces
we extend the definition of \eqref{deftrsp} to all values $m\leqslant 2n$;
see \cite[Sec.~3.3]{m:ff} for more details.

Recall that the factorial elementary symmetric polynomials
$e_r(l^2_1,\dots, l^2_n\tss|\tss a)$ with the
sequence $a$ defined as in \eqref{seqa} with $\ve=1$
are algebraically independent generators of the algebra
$\CC[l_1^2,\dots,l_n^2]^{\Sym_n}$.
Now we let $m$ be fixed and let $n$ run over integer values $\geqslant m/2$.
Set $k=\lfloor m/2\rfloor$. It is clear from the above argument
and from the explicit formulas for the symmetrizer $S^{(m)}$ that
the Harish-Chandra image
of the Casimir element \eqref{deftrsp} can be written
as the following linear combination
\ben
\sum_{\pi} c^{(m)}_{\pi}(u_1,\dots,u_m)\tss
\prod_{r=1}^k e_r(l^2_1,\dots, l^2_n\tss|\tss a)^{p_r},
\een
where $\pi$ runs over the $k$-tuples of nonnegative integers
$\pi=(p_1,\dots,p_k)$ satisfying the condition
$p_1+2\tss p_2+\dots+k\tss p_k\leqslant m/2$.
Moreover,
the $c^{(m)}_{\pi}(u_1,\dots,u_m)$
are symmetric polynomials in $u_1,\dots,u_m$ whose coefficients
are rational functions in $n$.
Since a rational function is uniquely determined by its values
at infinitely many points, the Harish-Chandra image of the Casimir element
\eqref{deftrsp} is uniquely determined by its values at 
infinitely many values of $n$.

\bth\label{thm:scaseven}
For any $1\leqslant k\leqslant n$ the image of the Casimir element
\beql{sgenbth}
\frac{n-k+1}{n-2k+1}\ts\tr\ts (F_1-k+1)
\dots (F_{2k}+k)\ts S^{(2k)}\in \U(\spa_{2n})^{\Spr_{2n}}
\eeq
under the Harish-Chandra isomorphism coincides with
$(-1)^k\ts e_k(l^2_1,\dots, l^2_n\tss|\tss a)$.
\eth

\bpf
As explained above,
it will be sufficient to prove the statement for any fixed $k$ under
the assumption that $n\geqslant 2k$. We will use the same argument
as in the proof of Theorem~\ref{thm:caseven}. Let $D_k$ denote
the Casimir element \eqref{sgenbth}. Relying on
Proposition~\ref{prop:charz} we show that $D_k$ vanishes
in all representations $V(\la)$ of $\Spr_{2n}$, where the partitions
$\la$ satisfy $|\la|<k$. Set $r=|\la|$
and consider the action of the Lie algebra $\spa_{2n}$
on $\CC^{2n}$ defined by
\ben
F_{ij}\mapsto -e_{ji}+\ve_i\tss\ve_j\tss e_{i'j'},\qquad 1\leqslant i,j\leqslant 2n.
\een
Under the corresponding representation in the tensor space
\ben
\vp:\U(\spa_{2n})\ot\End(\CC^{2n})^{\ot 2k}\to
\End(\CC^{2n})^{\ot r}\ot \End(\CC^{2n})^{\ot 2k}
\een
we have
\beql{simfa}
\vp(F_a)= \sum_{b=1}^r (-P_{\tss b\ts\ts r+a}+Q_{\tss b\ts\ts r+a}),\qquad a=1,\dots,2k,
\eeq
where we use the operators \eqref{pdef} and \eqref{qdefsp}.
The vanishing condition of Proposition~\ref{prop:charz} will be
verified if we show that if $|\la|<k$ then
\beql{setfze}
\big(\vp(F_1)-k+1\big)\dots \big(\vp(F_{2k})+k\big)\ts S^{\tss\prime(2k)}=0,
\eeq
where $S^{\tss\prime(2k)}$ denotes the image of the symmetrizer in the Brauer
algebra $\Bc_{2k}(-2n)$ acting of the last $2k$ copies of the tensor product space
$(\CC^{2n})^{\ot(r+2k)}$. Due to \eqref{simfa} and relations \eqref{mulse},
the desired identity \eqref{setfze} can be written in the form
\ben
(x_{r+1}-k+1)\dots (x_{r+2k}-k+1)\ts S^{\tss\prime(2k)}=0,
\een
where $x_{r+1},\dots,x_{r+2k}$ denote the images of the Jucys--Murphy
elements \eqref{jmdef} of the Brauer algebra $\Bc_{r+2k}(-2n)$
under its action \eqref{braactsp} in $(\CC^{2n})^{\ot(r+2k)}$.
It is verified by the same argument as
in the proof of Theorem~\ref{thm:caseven}.
The leading term of $\chi(D_k)$
coincides with $(-1)^k e_k(\la^2_1,\dots, \la^2_n)$ as shown in
\cite{m:ff}
and \cite{mr:cm}.
\epf

For $1\leqslant m\leqslant 2n$ consider the polynomials
in a variable $u$ whose coefficients are Casimir elements for $\spa_{2n}$ given by
\ben
D_{m}(u)=\frac{n-m/2+1}{n-m+1}\ts
\tr\tss \Big(F_1+u+\frac{m-1}{2}\Big)\Big(F_2+u+\frac{m-3}{2}\Big)
\dots \Big(F_{m}+u-\frac{m-1}{2}\Big)\ts S^{(m)}.
\een

\bco\label{cor:sdpol}
For the images
under the Harish-Chandra isomorphism we have
\ben
\chi:D_{m}(u)\mapsto
\ts\sum_{r=0}^{\lfloor\frac{m}{2}\rfloor}\ts(-1)^{r}\ts
\binom{2n-2r+1}{m-2r}\ts e_{r}(l^2_1,\dots, l^2_n\tss|\tss a)
\ts \prod_{i=0}^{m-2r-1}\Big(u-\frac{m-1}{2}+r+i\Big)
\een
and
\begin{multline}
\chi:\frac{n-m/2+1}{n-m+1}\ts\tr\ts
(-\di_t+F_1\tss t^{-1})\dots (-\di_t+F_m\tss t^{-1})\ts S^{(m)}
\mapsto \\
\sum_{r=0}^{\lfloor\frac{m}{2}\rfloor}\ts(-1)^r\ts
\binom{2n-2r+1}{m-2r}\ts e_{r}(l^2_1,\dots, l^2_n\tss|\tss a)
\ts t^{-2r}\ts (-\di_t+r\tss t^{-1})^{m-2r}.
\non
\end{multline}
\eco

\bpf
Exactly
as in the proof of Corollary~\ref{cor:dpol}, we use Lemma~\ref{lem:invpe}
to verify the relations $D_{m}(u)=(-1)^m\ts D_{m}(-u)$ and
\ben
D_m(u+1/2)-D_m(u-1/2)=(2n-m+2)\ts D_{m-1}(u),\qquad m\geqslant 1,
\een
with $D_0(u)=1$. Since
the same relation is satisfied by
the polynomials which are claimed to be the images $\chi(D_m(u))$,
the statement follows from Theorem~\ref{thm:scaseven}.
\epf

Finally, we obtain analogues of Theorem~\ref{thm:scaseven}
and Corollary~\ref{cor:sdpol}, where
$S^{(m)}$ is replaced by $A^{(m)}$.
We keep the
notation $A^{(m)}$ for the image
of the anti-symmetrizer under the action \eqref{braactsp}.
Note that due to the minus signs in that formula, the operator
$A^{(m)}$ acts in the tensor space \eqref{tenprk} as
the {\it symmetrization operator\/}.
The corresponding
Casimir elements coincide with those constructed
in \cite{i:tp} and \cite{mn:ci}.

\bth\label{thm:sascaseven}
For any $k\geqslant 1$ the image of the Casimir element
\beql{sasgenbth}
\tr\ts (F_1+k-1)\dots (F_{2k}-k)\ts A^{(2k)}\in \U(\spa_{2n})^{\Spr_{2n}}
\eeq
under the Harish-Chandra isomorphism coincides with
$h_k(l^2_1,\dots, l^2_n\tss|\tss a)$.
\eth

\bpf
Denote the Casimir element \eqref{sasgenbth} by $C_k$. We use
Proposition~\ref{prop:charz} and show that $C_k$ vanishes
in all representations $V(\la)$ of $\Spr_{2n}$, where the partitions
$\la$ satisfy $|\la|<k$. Using the same realization of $V(\la)$
as in the proof of Theorem~\ref{thm:scaseven}, we come to showing
the following analogue of \eqref{setfze}:
\ben
\big(\vp(F_1)+k-1\big)\dots \big(\vp(F_{2k})-k\big)\ts A^{\tss\prime(2k)}=0.
\een
We verify that
the image of the operator
$(-x_{r+1}+k-1)\dots (-x_{r+2k}+k-1)$
on the vector space $(\CC^{2n})^{\ot(r+2k)}$
is contained in a direct sum of representations $V(\nu)$ of $\Spr_{2n}$
with $\nu_1<k$,
while the image of the operator $A^{\tss\prime(2k)}$
is contained in a direct sum of representations $V(\nu)$ with $\nu_1\geqslant k+1$.
An easy calculation shows that the leading term of $\chi(C_k)$
coincides with $h_k(\la^2_1,\dots, \la^2_n)$.
\epf

Now consider the polynomials
in $u$ given by
\ben
C_{m}(u)=
\tr\tss \Big(F_1+u+\frac{m-1}{2}\Big)\Big(F_2+u+\frac{m-3}{2}\Big)
\dots \Big(F_{m}+u-\frac{m-1}{2}\Big)\ts A^{(m)}.
\een
Their coefficients are Casimir elements for $\spa_{2n}$.

\bco\label{cor:sasdpol}
For the images
under the Harish-Chandra isomorphism we have
\ben
\chi:C_{m}(u)\mapsto
\sum_{r=0}^{\lfloor\frac{m}{2}\rfloor}\ts
\binom{2n+m-1}{m-2r}\ts
h_{r}(l^2_1,\dots, l^2_n\tss|\tss a)
\ts \prod_{i=0}^{m-2r-1}\Big(u-\frac{m-1}{2}+r+i\Big)
\een
and
\begin{multline}
\chi:\tr\ts (-\di_t+F_1\tss t^{-1})
\dots (-\di_t+F_m\tss t^{-1})\ts A^{(m)}\mapsto\\
\sum_{r=0}^{\lfloor\frac{m}{2}\rfloor}\ts
\binom{2n+m-1}{m-2r}\ts
h_{r}(l^2_1,\dots, l^2_n\tss|\tss a)
\ts t^{-2r}\ts (-\di_t+r\tss t^{-1})^{m-2r}.
\non
\end{multline}
\eco

\bpf
Lemma~\ref{lem:invpe} implies
$C_{m}(u)=(-1)^m\ts C_{m}(-u)$ so that for the odd values $m=2k-1$ we have
$C_{2k-1}(0)=0$, while
$\chi(C_{2k}(1/2))=\chi(C_{2k}(-1/2))$ is found
from Theorem~\ref{thm:sascaseven}. Furthermore,
we verify that
\beql{sasdrecr}
C_m(u+1/2)-C_m(u-1/2)=(2n+m-1)\ts C_{m-1}(u),\qquad m\geqslant 1,
\eeq
where $C_0(u)=1$. This follows from Lemma~\ref{lem:invpe}
and the relation
for the partial trace of the operator $A^{(m)}$
found by
\ben
\tr_m\ts A^{(m)}=\frac{2n+m-1}{m}\ts
A^{(m-1)}
\een
thus verifying \eqref{sasdrecr}. The same relation is satisfied by
the polynomials which are claimed to be the images $\chi(C_m(u))$ as stated
in the corollary.
\epf

%\newpage


\begin{thebibliography}{99}


\bibitem{b:aw}
R. Brauer,
{\it On algebras which are connected with
the semisimple continuous groups},
Ann. Math. {\bf 38} (1937), 854--872.

\bibitem{d:ae}
{J. Dixmier},
{\it Alg\`ebres Enveloppantes},
{Gauthier-Villars, Paris},
1974.

\bibitem{hu:ci}
{R. Howe and T. Umeda},
{\it The Capelli identity, the double commutant theorem,
and multiplicity-free actions},
{Math. Ann.}
{\bf 290}
(1991),
569--619.

\bibitem{im:fp}
A. P. Isaev and A. I. Molev,
{\it Fusion procedure for the Brauer algebra},
St. Petersburg Math. J. {\bf 22} (2011), 437--446.

\bibitem{i:td}
{M. Itoh},
{\it Two determinants in the universal enveloping algebras
of the orthogonal Lie algebras}
J. Algebra {\bf 314} (2007), 479--506.

\bibitem{i:tp}
{M. Itoh},
{\it Two permanents in the universal enveloping algebras of
the symplectic Lie algebras}, Internat. J. Math. {\bf 20} (2009), 339--368.

%%\bibitem{k:dt}
%%{K. Koike}
%%{\it On the decomposition of tensor products of the representations of
%%the classical groups: by means of the universal characters},
%%Adv. Math. {\bf 74} (1989), 57--86.

\bibitem{lr:rh}
R. Leduc and A. Ram,
{\it A ribbon Hopf algebra approach to the
irreducible representations of
centralizer algebras: The Brauer, Birman-Wenzl and
type A Iwahori-Hecke algebras}, Adv. Math.
{\bf 125} (1997), 1--94.

\bibitem{m:sft}
{I. G. Macdonald},
{\it Schur functions: theme and variations}, in:
``Actes 28-e S\'eminaire Lotharingien", pp. 5--39,
{Publ. I.R.M.A. Strasbourg,\/} 1992, 498/S--27.

\bibitem{m:ff}
A. I. Molev,
{\it Feigin--Frenkel center in types $B$, $C$ and $D$},
Invent. math. (2012); {\tt doi: 10.1007/s00222-012-0390-7}.

\bibitem{mn:ci}
{A. Molev and M. Nazarov},
{\it Capelli identities for classical Lie algebras},
{Math. Ann.} {\bf 313} (1999), 315--357.

\bibitem{mr:cm}
A. I. Molev and N. Rozhkovskaya,
{\it Characteristic maps for the Brauer algebra},
{\tt arXiv:1112.0620}.

\bibitem{n:yo}
{M. Nazarov},
{\it Young's orthogonal form for Brauer's
centralizer algebra},
J. Algebra {\bf 182} (1996), 664--693.

\bibitem{n:ce}
{M. Nazarov},
{\it Capelli elements in the classical universal enveloping algebras}, in:
``Combinatorial methods in representation theory (Kyoto, 1998)'',
Adv. Stud. Pure Math. {\bf 28}, Kinokuniya, Tokyo, 2000,
pp. 261--285.

\bibitem{o:qi}
{A. Okounkov},
{\it Quantum immanants and higher Capelli identities},
{Transform. Groups} {\bf 1} (1996), 99--126.

\bibitem{oo:ss}
A. Okounkov and G. Olshanski,
{\it Shifted Schur functions},
{St.\,Petersburg Math. J.} {\bf 9} (1998), 239--300.

\bibitem{oo:ss2}
A. Okounkov and G. Olshanski,
{\it Shifted Schur functions II.
Binomial formula for characters
of classical groups and applications}, in:
``Kirillov's Seminar on Representation Theory''
(G.~Olshanski, Ed.),
{Amer. Math. Soc. Transl.}
{\bf 181},
AMS,
Providence, RI,
1998, pp. 245--271.

\bibitem{o:gs}
{G. Olshanski},
{\it Generalized symmetrization in enveloping algebras},
Transform. Groups {\bf 2} (1997), 197--213.

\bibitem{w:ce}
A. Wachi,
{\it Central elements in the universal enveloping algebras
for the split realization of the orthogonal Lie algebras},
Lett. Math. Phys. {\bf 77} (2006), 155--168.


\end{thebibliography}
\end{document}